\documentclass{article}
\usepackage{blindtext}
\usepackage[a4paper, total={6in, 8in}]{geometry}

\usepackage{graphicx}% Include figure files
\usepackage{subfig}
    \graphicspath{{./figures/}}
\usepackage{dcolumn}% Align table columns on decimal point
\usepackage{bm}% bold math
\usepackage[mathlines]{lineno}% Enable numbering of text and display math
\usepackage{tabularx,array,booktabs}
\usepackage{caption}
\usepackage{mathtools}
\usepackage{amsmath,amssymb}
\usepackage{soul}
\usepackage{lipsum}
\usepackage{varioref}
\usepackage{algorithm,algpseudocode}
\usepackage{ulem}
\usepackage{siunitx}
\usepackage{comment}
\usepackage{upgreek}
\usepackage{url}
    
\usepackage{calc}
\newcommand{\bb}[1]{\mathbf{#1}}
\newcommand{\bs}[1]{\boldsymbol{#1}}
\newcommand{\R}{\mathbb R}
\newcommand{\C}{\mathbb C}
\newcommand{\LL}{\mathcal L}

\newcommand{\Comsol}{COMSOL Multiphysics\textsuperscript{\textregistered}}

\begin{document}

\title{Optimal Design of Broadband, Low-Directivity Graded Index Acoustic Lenses for Underwater Applications}
\date{}

\author{Sebastiano Cominelli \thanks{Sebastiano Cominelli is a PhD Candidate at Politecnico di Milano, Department of Mechanical Engineering, Milano 20133, Italy (e-mail: sebastiano.cominelli@polimi.it), Corresponding author.}
\and
Francesco Braghin\thanks{Prof. Francesco Braghin is Full Professor at Politecnico di Milano, Department of Mechanical Engineering, Milano 20133, Italy (e-mail: francesco.braghin@polimi.it). }
}

\maketitle
\begin{abstract}
    Manipulating underwater pressure waves is crucial for marine exploration, as electromagnetic signals are strongly absorbed in water. However, the multi-path phenomenon complicates the accurate capture of acoustic waves by receivers. Although graded index lenses, based on metamaterials with smoothly varying properties, successfully focus pressure waves, they tend to have high directivity, which hinders practical application.

This work introduces three 2D acoustic lenses made from a metamaterial composed of solid inclusions in water. We propose an optimization scheme where the pressure dynamics is governed by Helmholtz's equation, with control parameters affecting each lens cell's density and bulk modulus. Through an appropriate cost function, the optimization encourages a broadband, low-directivity lens.\\
The large-scale optimization is solved using the Lagrangian approach, which provides an analytical expression for the cost gradient. This scheme avoids the need for a separate discretization step, allowing the design to transition directly from the desired smooth refractive index to a practical lattice structure. As a result, the optimized lens closely aligns with real-world behavior.

The homogenized numerical model is validated against finite elements, which considers acoustic-elastic coupling at the microstructure level. When homogenization holds, this approach proves to be an effective design tool for achieving broadband, low-directivity acoustic lenses.
\end{abstract}

% \begin{fmtext}
% %%%% Insert A head here

% \end{fmtext}

%\end{frontmatter}
\maketitle

\section{Introduction}
\label{sec:intro}
% \seba{Penso di sottometterlo a JASA, JSV o PRSA, cosa ne dici?}

\textit{We know more about the moon than our own oceans}~\cite{snelgrove2010discoveries}.
Even though this commonly made claim is controversial, it emphasizes how tough it is for mankind to explore a marine environment. Among several technical difficulties compared to the mainland, long-distance communications cannot rely on electromagnetic waves because water is a highly absorptive medium, causing significant attenuation and loss of signals strength~\cite{jiang2011electromagnetic}.
Conversely, sound travels efficiently in water because it experiences minimal attenuation
% and it penetrates various underwater obstacles, such as sediment and marine life
\cite{brown2003ray}. 

Pressure waves have disadvantages with respect to electromagnetic waves. For instance a limited data rate and high power consumption~\cite{zia2021state}. In addition, the pressure intensity progressively drops as the wave spreads through the ocean thus long-distance signals inherently suffer from low strength.

In this context, researchers from all over the world are developing innovative components of the transmission line, e.g., projectors, hydrophones, and modems~\cite{furqan2019underwater,jiang2016convert}. Among them, a lot of effort is spent on the receiver side to improve the signal-to-noise ratio, that is obtained through intensification of the pressure waves through acoustic lenses. The interest shown during the last two decades proves that this topic is far from being exhaustively explored~\cite{allam20203d,lin2009gradient,climente2010sound,martin2010sonic,zigoneanu2011design,li2014three,huang2021frequency,ma2022underwater,brambilla2024high}.

As it happens to light rays, an acoustic lens has the ability to focus pressure waves in a limited focal region. The performance of a lens strongly depends on the shape and on the properties of the device itself wherein the signal path is drawn by Snell's law.

Historically, several analytical expressions describing the refraction index of a lens have been derived. For instance, the fish-eye from Maxwell~\cite{whewell1854cambridge}, the L\"uneburg lens~\cite{luneburg1966mathematical}, and the hyperbolic secant profile are among the most famous distributions. Where the former two belong to classical optics and the latter has been more recently proved to reduce the aberration~\cite{lin2009gradient,su2017broadband}.
\\
Such devices, relying on smooth refractive index profiles, are commonly called gradient index (GRIN) lenses and are attained through gratings having slow space variation of the unit cell. In this regard, many experiments have been conducted on metamaterials made of periodic scatterers~\cite{zigoneanu2011design,su2017broadband,huang2021frequency}.

% - Some literature for acoustics communication: describing the innovation of almost each paper
%~\cite{ma2022underwater} Planar lens for frequencies between 7 and 10 KHz
%~\cite{su2017broadband} experimental two dimensional lens realizing the hyperbolic secant profile
% \seba{Norris in~\cite{gokhale2012special} ipotizza un concentrator, ma quello non serve al nostro scopo perché non cambia l'ampiezza dell'onda}

% - general introduction on Metamaterials
% Gradient-index (GRIN)
% - Some examples of metamats for focusing and communicate (superlenses, orbital angular momentum)

All the acoustic lenses proposed in the literature have been designed to operate effectively in controlled environments, where the sound source and the lens are precisely located in space.
This is especially the case of the recent work by Ma~\textit{et al.}~\cite{ma2022underwater}, where experiments are conducted on a thin planar lens made of Helmholtz's resonators: the signal passing through the device experiences a precise phase delay such that it concentrates in the focal point.

However, such settings cannot apply to a marine environment, where the so-called multi-path phenomenon makes the direction of propagation hard to know \textit{a priori}~\cite{brown2003ray}. Marine water shows sensible temperature gradients that, along with solid inclusions, perturb the propagation of pressure waves causing the signal to follow diverse routes and reflections before reaching the lens.
Even if such a phenomenon is negligible where only near field effects are accounted, the sound propagation at the kilometer scale cannot disregard it. Moreover, in long distance communication the piece of information related to reciprocal position between sender and receiver is usually not readily available, hence even the most advanced prediction on ray propagation is pointless.

To this end, a design based on space coordinate transform becomes appealing since it leads to broadband and omni-directional performance~\cite{pendry2000negative,rahm2008design,gokhale2012special}. This mathematical tool recasts a chosen space transformation as peculiar properties of the device. However, the behavior of an anisotropic fluid is often required and, in case of a space folding, negative refraction index materials (NIM) are needed~\cite{veselago1967electrodynamics,pendry2000negative,imamura2004negative,cominelli2024isospectral}. Some materials exploiting anisotropy or negative effective properties have been tested~\cite{antonakakis2014gratings,shelby2001experimental}, but they are not practical.
Then, the serious difficulties in attaining the anisotropy required make such devices out of reach.
On the other hand, large scale optimization techniques have been proved to be effective for controlling waves field in a chosen frequency band and angular aperture~\cite{cominelli2022design,cominelli2023optimal,chen2021optimal}.

% \cbox{dire da qualche parte che è lente/risuonatore}
In this work we recast the design of acoustic lenses as an Optimal Control Problem (OCP) where the control parameters are the space varying physical properties of the lens and the cost is a measure of the pressure gain achieved in the lens focus (Section~\ref{sec:physics}). The OCP is constraint by an elliptical Partial Differential Equation (PDE), namely the acoustic Helmholtz equation for inhomogeneous, isotropic fluids. \\
We choose a 2D lattice available in literature \cite{cominelli2022design} to build the lens. Thus a piece-wise constant control space is introduced to account for its discrete nature and further constraints are imposed to limit the control to the set of attainable properties (Section~\ref{sec:sizing}). 
Given the large scale nature of the problem, we derive first order necessary conditions by using a Lagrangian approach and the numerical problem is solved relying on a gradient-based iterative method.

Initially, we present the optimized sizing of a directional lens to demonstrate the attainable focusing effect based on the lens size and operating frequencies.\\
In Section \ref{sec:multi-directional}, we address a robust Optimal Control Problem (OCP) to design a broad-frequency, low-directivity device. Three lenses are optimized for different central frequencies, with a focus on minimizing directivity to mitigate issues arising from multipath phenomenon. We analyze the results by examining the transfer function between the incident wave and the pressure signal detected at the lens foci.
Through a deeper analysis of the microstructure response, the three lenses are tested against the low frequency homogenization limits.\\
Finally, conclusive remarks are proposed and the future perspectives of this work are discussed in~Section \ref{sec:conclusion}.

% Note that a directional device allows to radiate a concentrated and focused beam if it is applied to a source antenna. Such a dual configuration finds its application in the case a preferred direction is given.

\section{Problem Statement}
\label{sec:physics}

An ideal lens capable to improve acoustic sensing should focus pressure signals coming from a wide range of directions and characterized by a broadband frequency content. This two requirements lend themselves to the choice of a circular annular lens able to receive pressure waves all around and concentrate them in its centre, where an acoustic sensor is hosted.
% Note that in underwater applications the most used pressure transducers are hydrophones, a family of sensors that rely on piezoelectric materials to measure the pressure field; for that,
In the following, the sensor will be equivalently referred to as hydrophone or transducer.

\begin{figure}
	\centering
	\includegraphics[width=0.35\textwidth,trim=0 15 0 0]{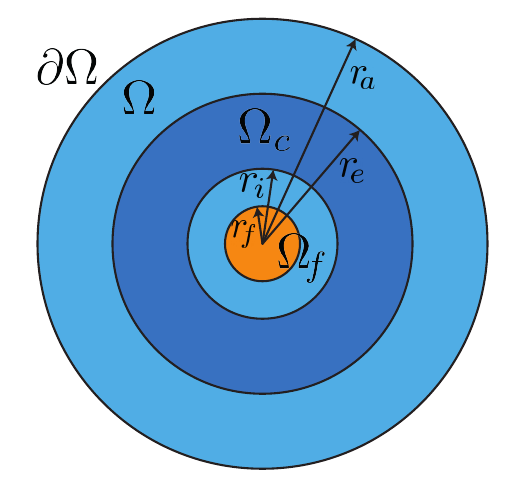}
	\caption[Domain definition]{Schematic representation of the computational domain $\Omega$. The lens occupies the control region $\Omega_c$, filled by inhomogeneous properties accounting for the grating. The background fluid occupies the remaining domains. $\Omega_f$ is the focal region where the sensor is placed.}
	\label{fig:domain definition}
\end{figure}

With reference to Figure~\ref{fig:domain definition}, two subdomains are recognised in the computational domain $\Omega \subset \R^2$: $\Omega_c$ is the control domain occupied by the lens and it is filled with an inhomogeneous yet isotropic medium able to represent the grating; $\Omega_f$ is the portion of $\Omega$ where we aim at maximizing the pressure intensity, namely the focal region hosting a transducer. 

In order to have a complete model, one should take into account the dynamics of a hydrophone and its coupling with water. However, a good sensing is guaranteed if the transducer size is much smaller that the signal wavelength, such that it experiences an in-phase pressure on its surface. This allows us to disregard its confined loading effect for the sake of the model simplicity. Also, this leads to results comparable with the existing literature, where only the undisturbed field is analysed.
Inherently, the size of $\Omega_f$ is chosen slightly larger than that of the hydrophone. This enhances the sensing and allows to choose the best hydrophone position \textit{a posteriori}.

More than that, the lens accounts for a clearance between the focal region $\Omega_f$ and the lens itself $\Omega_c$. This volume and the surrounding ambient are filled with water. The model domain is bounded by $\partial\Omega$ for computational reasons only and the free field propagation is approximated by artificial absorbing conditions.
We consider water having density $\rho_0 = \SI{998}{\kilo\gram\per\cubic\meter}$ and bulk modulus $\kappa_0=\SI{2.2}{\mega\pascal}$, where $\rho_0$ and $\kappa_0$ are also referred to as the reference mass and stiffness. On the other hand, the physical properties in the domain $\Omega_c$ are addressed as control functions and denoted by $\rho(\bb x)$ and $\kappa(\bb x)$, functions of the space coordinates $\bb x \in\Omega_c$. 

When the system is forced by a time harmonic wave, the steady-state acoustic pressure $P(\bb{x},t)$ is 
governed by
% formulated as $\Re\{p(\bb{x})\,e^{ j \omega t}\}$, where $\Re\{\cdot\}$ denotes the real part of its argument. The complex amplitude $p(\bb{x})$ satisfies
the Helmholtz equation for inhomogeneous media \cite{bergmann1946wave}:
\begin{equation}	\label{eq:Helmholtz frequency}
	\nabla\cdot\big(a(\bb{x}) \nabla p(\bb{x})\big) = -b(\bb{x})\,\omega^2 p(\bb{x}),
\end{equation}
where $p(\bb{x}) \in \C$ is the pressure field phasor,  $P(\bb x,t)=\Re\{p(\bb{x})\,e^{ j \omega t}\}$, and $\omega$ the circular frequency of the forcing wave. The coefficients $a$ and $b$ are defined as $a \coloneqq \rho^{-1}$ and $b\coloneqq\kappa^{-1}$, where $\rho$ is the local mass density and $\kappa$ the local bulk modulus. The definition of $a$ and $b$ will turn out to be useful in manipulating \eqref{eq:Helmholtz frequency} and setting up the resulting OCP. \\
Notice that \eqref{eq:Helmholtz frequency} governs the dynamics of an ideal, inviscid, still fluid at homogeneous constant temperature. Such an acoustic wave equation is not able to capture typical phenomena of ocean acoustics (absorption, multipath propagation due to temperature gradients, scattering by air bubbles\dots \cite{brown2003ray}), yet it is a valid choice to model the near-field propagation around the lens---at the end of the long acoustic path through the ocean.
% Moreover, it is commonly adopted in the literature related to underwater acoustic lenses \cite{allam20203d,lin2009gradient,climente2010sound} because it is sufficiently accurate to allow for the design of experimental demonstrators.
\\
% As we anticipated in the introduction, we target the design of a lens which is robust to the incident direction, thus a more accurate model would be superfluous.\\
We also note that the acoustic wave equation has the same analytical structure of the equations describing out of shear horizontal (SH) elastic waves in solids and transverse electric (TE) and transverse magnetic (TM) modes in electromagnetism \cite{banerjee2011introduction}. It follows that, upon suitable reinterpretation of \eqref{eq:Helmholtz frequency}, the following results can be readily related to such types of incident radiation.

The total pressure field can be decomposed into an incident and a scattered field, $p(\bb x)=p^s(\bb x)+p^i(\bb x)$, where $p^i$ is the solution of the Helmholtz equation obtained considering a homogeneous fluid without obstacles, also known as background field. Considering that pressure waves travel for long distances in oceans before reaching the lens, we exploit the far field approximation and consider all incident waves as plane, i.e., $p^i(\bb x) = e^{-j \bb k\cdot\bb x}$, where $\bb k$ the wave vector defining the direction of the incident wave, $|\bb k| = k = \omega/c_{0}$ is the homogeneous wave number, $c_0=\sqrt{\kappa_0/\rho_0}$ is the undisturbed sound velocity in water.

Following \cite{cominelli2022design}, we define the OCP
\begin{align}	\label{eq:cost functional}
	\min_{v,u,p^s} J(v,u,p^s) &=
        \frac{\sigma}{2} \left|(v,u)\right|^2_{\mathcal U}
        % \frac{\lambda_v}{2} \int_{\Omega_c} \left|v\right|^2_{\mathcal U} \,d\Omega
        % + \frac{\lambda_u}{2} \int_{\Omega_c} \left|u\right|^2_{\mathcal U} \,d\Omega
        - \frac{1}{2}G(p^s)
        % + \frac{1}{2} \int_{\Omega_a}{\bar{p}_{s} p^s	\,d\Omega}
	\\
	s.t. &
	\begin{dcases} 	\label{eq:state_dyn}
		-\nabla \cdot (a\nabla p^s) - b \omega^2 p^s = f	& \text{in } \Omega
		\\
	    a\nabla p^s\cdot\mathbf{n} + \alpha p^s  =  0 	& \text{on } \partial\Omega
	\end{dcases}
\end{align}
where
\begin{equation}
    \begin{dcases}
		f = \omega^2(b-b_0)p^i +\nabla\cdot[(a-a_0)\nabla p^i]
		\\
		\alpha = a\Big(jk + \frac{1}{2r_a}\Big)
    \end{dcases}
\end{equation}
and the perturbed material properties depend on the control functions $u(\bb x)$ and $v(\bb x)$ such that:
\begin{align}
    a &= a_0\,e^{-v}, & b &= b_0\,e^{-u}.
\end{align}
The positive constant $\sigma$ increases the problem convexity limiting the control effort.

In view of setting a well posed OPC, suitable functional spaces for state and control functions are defined: we select the complex-valued Hilbert space $H^1(\Omega,\C)$ as the state space, that is $\mathcal{V}=H^1(\Omega,\C)$. The boundary value problem \eqref{eq:state_dyn} is well-posed as long as its coefficients $\kappa$ and $\rho$ are bounded and positive \cite{colton1998inverse}. The exponential modulation we choose for the background properties allows us to select as control space $\mathcal U=L^{\infty}(\Omega_c,\R)^2\cap H^1(\Omega_c,\R)$, that is the space of real-valued two-dimensional vector functions which are essentially bounded and have a bounded gradient. In other words, for each $\bb{x} \in \Omega_c$ we associate a real-valued control pair $(v(\bb x),u(\bb x))$ whose elements are bounded and with derivatives bounded in every portion of $\Omega_c$ with a non-null measure.
This condition permits the norm $\left|v,u\right|_{\mathcal U}$ to weight the space fluctuations of the properties, i.e.,
\begin{equation}
    \left|v,u\right|^2_{\mathcal U}\coloneqq \int_{\Omega_c} (v^2+\left|\nabla v\right|^2 + u^2+\left|\nabla u\right|^2 )\,d\Omega.
\end{equation}
Such a norm enforces the control to have a further regularity, useful in view of the lens realization: by definition, a GRIN lens is made up of a metamaterial composed by cells whose properties smoothly vary in space. The effective property of each cell is computed under the assumption of an infinite periodic lattice, that is far from the real case, but to contain the properties gradient makes the lens closer to such hypothesis.

Note that the robin boundary condition imposed on $\partial\Omega$ is the first-order approximation of the well-known Sommerfeld radiation condition for 2D domains \cite{schot1992eighty,bayliss1980radiation}, where $r_a$ is the radius of $\Gamma_e$.\\
This approximation guarantees a perfect absorption when the wave direction is normal to the surface $\partial\Omega$, and a reflection coefficient smaller than \SI{1}{\percent} for waves impinging within $\SI{30}{\degree}$ with respect to the boundary normal \cite{Comsolblog}.
% reflection coefficient that increases up to 0.01 when the angle between the propagation direction and the normal is $\SI{30}{degree}$ and to 0.1 when it is $60^o$ \cite{Comsolblog}.
We assume that the lens is close to the center of $\Omega$ such that any scattered wave has a nearly orthogonal incidence on the external surface, resulting in contained artificial reflections.
% This allows to consider such an approximation a computationally efficient implementation when compared to more sophisticated domain truncations \cite{negri2015efficient,bangtsson2003shape}.

The functional $G\colon \mathcal V\to\R$ measures the integral intensity gain as
\begin{equation}
    G(p^s) \coloneqq \int_{\Omega_f} (\bar p^s +\bar p^i)(p^s+p^i) \,d\Omega .
\end{equation}
Note that the scalar value given by $G$ is proportional to the mean pressure intensity in $\Omega_f$ and targets the lens performance such that a high intensity in the focal domain $\Omega_f$ reduces the cost $J$.
$\bar p^s$ and $\bar p^i$ represent the complex conjugate of $p^s$ and $p^i$, respectively.

We now derive a set of first-order optimality conditions applying Lagrange's method~\cite{Trol2010}, such that the explicit expression for the gradient of the cost functional in the continuous setting is computed.
The Lagrangian $\mathcal{L}\colon \mathcal{V} \times \mathcal{U} \times \mathcal{W}^{*} \to \R$ is formed as:
\begin{equation}\label{eq:Last L}
	\LL \coloneqq J
	+ \Re\Big\{\int_{\Omega} (\nabla \cdot (a\nabla p^s) + b \omega^2 p^s + f)\bar\lambda \, d\Omega 
    -\int_{\partial\Omega} (a\nabla p^s\cdot\bb n+\alpha p^s)\bar\lambda\,d\partial\Omega
    \Big\}
\end{equation}
where the adjoint function $\lambda\colon\Omega \to \C$ belongs to $H^1(\Omega,\C)$, i.e., we can identify $\mathcal{W}^{*} = H^{1}(\Omega,\C)$.

% \noindent In the following we will omit the explicit dependence on the space variable when it is clear from the context.

The first order optimality conditions are computed by requiring that about the optimal solution the Lagrangian has to be stationary with respect to small variations of the state, the adjoint and the control functions. Thus, computing G\^ateaux derivative of \eqref{eq:Last L} with respect to the state $p^s$ and equating it to zero, one obtains the following strong formulation:
\begin{equation}\label{eq:PDE adjoint}
    \begin{dcases}
        -\nabla \cdot (a\nabla \lambda) - b\omega^2 \lambda = -(p^s+p^i)\, \chi_{\Omega_f}		& \text{in } \Omega \\
        \nabla \lambda\cdot\mathbf{n} + \bar\alpha \lambda  =  0 		& \text{on } \partial\Omega
    \end{dcases},
\end{equation}
which is the adjoint dynamics, being $\chi_{\Omega_f}(\bb{x})$ the indicator function of the domain $\Omega_f$. Then, by computing the variation with respect to the controls $v$ and $u$, the gradients of the reduced cost are obtained
\begin{equation}\label{red_grad}
\begin{aligned}
    \nabla J_v &= \sigma (v+\Delta v) + \Re\big\{a\nabla(p^s + p^i)\cdot\nabla \bar\lambda\big\}=0 \\
    \nabla J_u &= \sigma (u+\Delta u) -  \Re\big\{b\omega^2 (p^s + p^i) \bar\lambda\big\}=0,
\end{aligned}
\quad\text{in }\Omega_c
\end{equation}
where the natural boundary conditions $\nabla v\cdot\bb n = 0$ and $\nabla u\cdot\bb n =0$ on $\partial\Omega_c$ are assumed on the control functions.\\
For a detailed derivation, the reader is referred to~\cite{cominelli2022design}, but particular attention is posed here on the right hand side of the state equation \eqref{eq:PDE adjoint}. It derives from the first variation of the functional $G(p^s)$, that can be computed as follows
\begin{equation}
    D_{p^s} G\,[\varphi] = \int_{\Omega_f}(\bar p^s+\bar p^i) \,\varphi\,d\Omega
    =0\quad\forall\varphi\in\mathcal V.
\end{equation}
In the following, a different definition of $J$ will change the adjoint equation.

\section{Sizing}
\label{sec:sizing}
The size of the lens has a significant impact on its performance for two main reasons.
First, a larger device interacts with a broader area of the signal, potentially allowing for greater gain at the focal point. This makes it possible to increase performance, but at the cost of a bulkier, less manageable lens. A common guideline is that its diameter should be at least as large as the signal wavelength; otherwise, the pressure waves does not effectively interact with the device. 
Second, a larger lens provides more flexibility in design, offering additional degrees of freedom because the control domain is the entire device.\\
Anyway, the refractive index is limited by the available materials, affecting the possible curvature of the acoustic path.
Thus, balancing lens size with these constraints is a critical initial step in design.

\begin{figure}
    \centering
    \subfloat[]{\includegraphics[width=0.62\textwidth,trim=0 0 0 0]{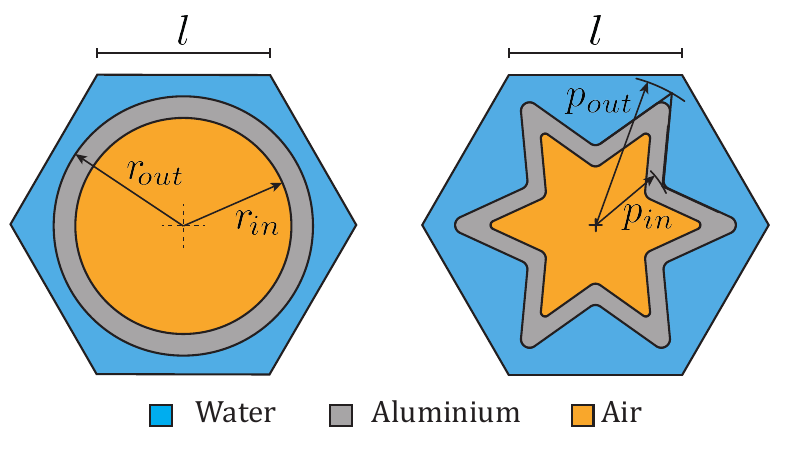}\label{fig:metamat-a}}
    % \subfloat[]{\includegraphics[width=0.36\textwidth,trim=0 -5 0 0,clip]{Figures/FrontCom 5.pdf}\label{fig:metamat-b}}
    \quad
    \subfloat[]{\includegraphics[width=0.28\textwidth,trim=0 0 -20 0,clip]{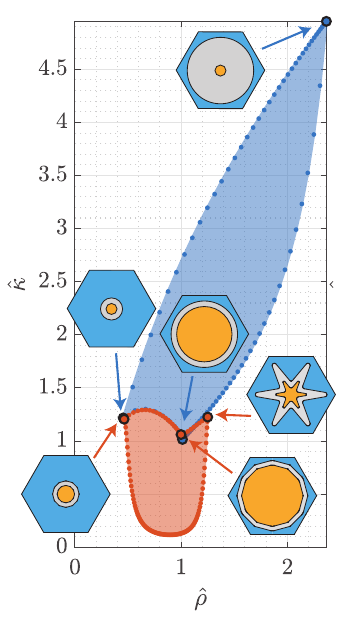}\label{fig:metamat-c}}
    \caption{ \protect\subref{fig:metamat-a}
    % and \protect\subref{fig:metamat-b}
    schematic of the two parametric families of cells composing the grating; \protect\subref{fig:metamat-c} the properties attainable by varying the parameters of this two cells.  
    The union of the blue and the orange areas outline the constraint for the control space for the optimization problem. Please note that the star shaped cell has 12 tips even if the drawings show 6 tips for the sake of simplicity.
    }
    \label{fig:metamat}
\end{figure}

The work here presented is based on some results by Cominelli~\textit{et al.}~\cite{cominelli2022design}. In particular, the metamaterial family analyzed in that work is adopted to reproduce an isotropic inhomogeneous media with which we can control the pressure field within the computational domain.
Considering the scale of the metamaterial small enough with respect to the lens size, we can leverage the long-wavelength homogenization and approximate the local dynamics of the microstructure considering only the quasi-static properties of each cell \cite{laude2015phononic}.

With reference to Figure~\ref{fig:metamat-a}, the cell edge is assumed $l=\SI{5}{\milli\meter}$ and the following geometric constraints hold:
\begin{equation}
    \begin{aligned}
        r_{in}&\geq l/5, &
        r_{out}&\leq \big(\sqrt{3}/2-1/5\big)\,l, &
        r_{out}&\geq r_{in} + l/20,
        \\
        p_{in} &\geq l/5, &
        p_{out} &\leq \big(\sqrt{3}/2-1/5\big)\,l, &
        p_{out} &\geq p_{in}.
    \end{aligned}
\end{equation}
The fillet radii and the thickness of the wall of the star shaped inclusion are chosen to be $l/40$ and $l/20$, respectively.\\
The attainable properties are outlined through the exploration of all possible parameters, and are shown in Figure~\ref{fig:metamat-c} in the plane $\hat\rho\times\hat\kappa$ as a blue and an orange regions, corresponding to the circular crown and to the star shaped cells, respectively. Note that $\hat\rho\coloneqq\rho/\rho_0$ and $\hat\kappa\coloneqq\kappa/\kappa_0$ are the density and the bulk modulus normalized with respect to the water properties.

The control space $\mathcal U$ is limited to the control pairs $(v(\bb x),u(\bb x))$ such that the property pairs $(e^v,e^u)$ belong to the union of these two colored regions on the plane $\hat\rho\times\hat\kappa$.

The dispersion relation of the cells show a linear relation when the reduced frequency is such that $\hat f=f\,l/c_0 < 0.09$ \cite{cominelli2022design}, i.e., $f<\SI{25}{\kilo\hertz}$ for $l=\SI{5}{\milli\meter}$.
However, such a theoretical limit holds for an infinite and periodic lattice while we target a finite size device whose microstructure is a smoothly varying grating.
Thus, such a limit is tested against three frequency bands centered on $f_L=\qty{10}{\kilo\hertz}$, $f_M=\qty{25}{\kilo\hertz}$, and $f_H=\qty{40}{\kilo\hertz}$ referred to as low, mid and high frequency, respectively. Considering the sound speed in water equal to $c_0=\SI{1485}{\meter\per\second}$, the corresponding wavelengths are $\lambda_L=\qty{148.5}{\milli\meter}$, $\lambda_M=\qty{59.4}{\milli\meter}$, and $\lambda_H=\qty{37.1}{\milli\meter}$.
% \qtyrange{9}{11}{\kilo\hertz}, \qtyrange{24}{36}{\kilo\hertz}, and \qtyrange{39}{41}{\kilo\hertz}. Note that the absolute bandwidth is kept constant, but not

The characteristic size of a typical omni-directional hydrophone intended for such a frequency range is about \qty{10}{mm}, see e.g.\ \cite{benthowave}.
Thus, we set $r_f=\qty{10}{\milli\meter}$ and $r_i=\qty{20}{\milli\meter}$. This choice allows the device to host a transducer in its center with a clearance, ensuring the sensing surface to be immersed within the focal region.
For the sake of simplicity, the lengths $l$, $r_f$, and $r_i$ are unchanged throughout this work.

A precise control on the pressure field is obtained through a fine tuning of the lens properties. Thus, the cell size plays a relevant role, even if it is negligible when $l/\lambda$ is small. An ideal case is when the properties of the lens are allowed to vary smoothly in space within the admissible control set $\mathcal U$, that corresponds to the limit $l/\lambda\to0$.

When $r_f,\,r_i\ll\lambda$ for any considered frequency, the design problem mainly depends on the geometrical ratio $r_e/\lambda$, relating the lens size to the wavelength.
So, if we completely disregard the dependence of the problem on $r_f$ and $r_i$ and consider, for instance, the low frequency only, the performance attainable for a generic frequency is readily obtained by looking at the results corresponding to the desired ratio $r_e/\lambda_L$.

Given the former assumptions, the OCP \eqref{eq:cost functional}-\eqref{eq:state_dyn} is discretized using second order Finite Elements (FEs) supported by a triangular mesh with characteristic size $\lambda/10$. The same discretization is adopted for the state, the adjoint and the control fields.
An iterative steepest descend method based on Armijo backtracking line-search is adopted and the constraints on the control space $\mathcal U$ are imposed using a standard Projected Gradient method. The algorithm is stopped when the length of one step or the module of the gradient becomes smaller then a fixed tolerance. The reader is referred to \cite{nocedal2006numerical,cominelli2022design} for a more detailed description.\\
Let us consider the low frequency $f_L$. The constrained problem is solved for different choices of the external radius $r_e$, such that $r_e\in[\qty{40}{\milli\meter},\qty{600}{\milli\meter}]$, corresponding to $r_e/\lambda_L\in[0.25,4]$. This parametric sweep gives a rough estimation of the achievable gain granted by the lens. The weight on the control effort is chosen to be equal to $\sigma=1$ to improve the stability of the solution limiting rapid variations of the properties. Indeed, the lens could take a huge benefit from them, but their realization is tough.

\begin{figure}
    \centering
    % \subfloat[]{\includegraphics[height=0.33\textwidth,trim=45 0 20 0,clip]{Figures/Analisi dimensione lenti 172 kHz 03-May-2024_p_mean.pdf}\label{fig:sizing costs-a}}
    \subfloat[]{\includegraphics[width=0.33\textwidth,trim=45 0 55 0,clip]{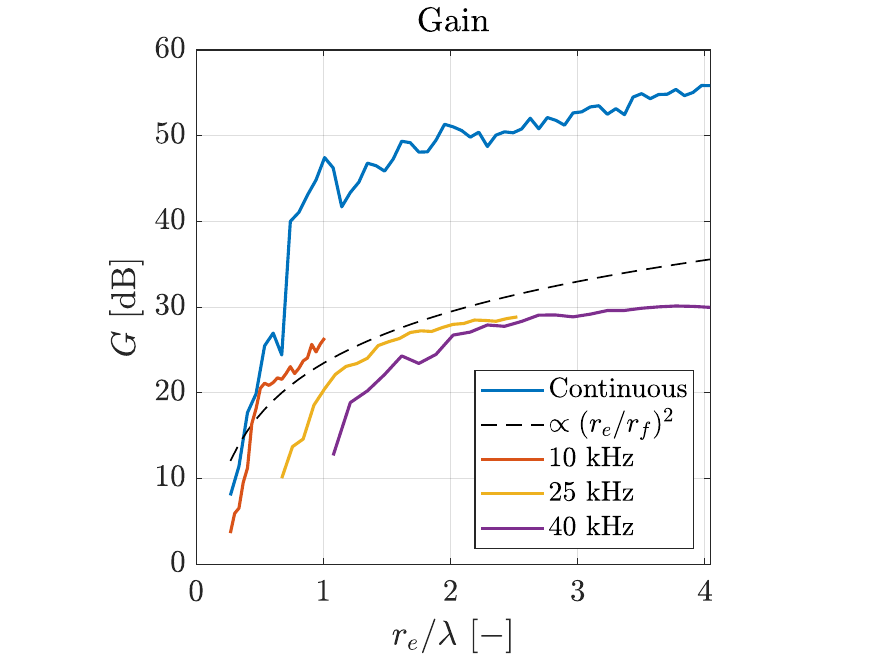}\label{fig:sizing costs-a}}
    \subfloat[]{\includegraphics[height=0.33\textwidth,trim=45 0 55 0,clip]{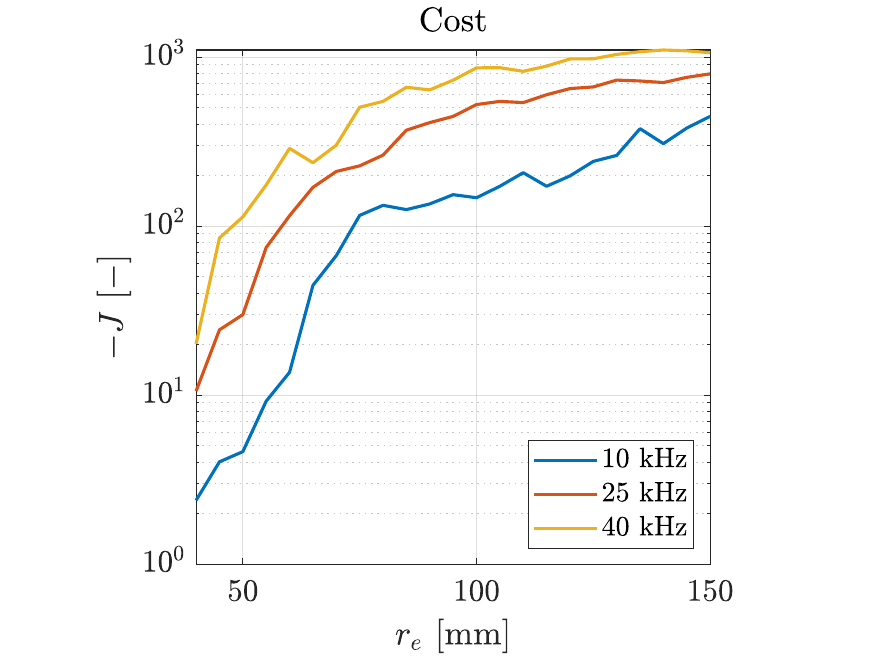}\label{fig:sizing costs-b}}
    \subfloat[]{\includegraphics[height=0.33\textwidth,trim=45 0 55 0,clip]{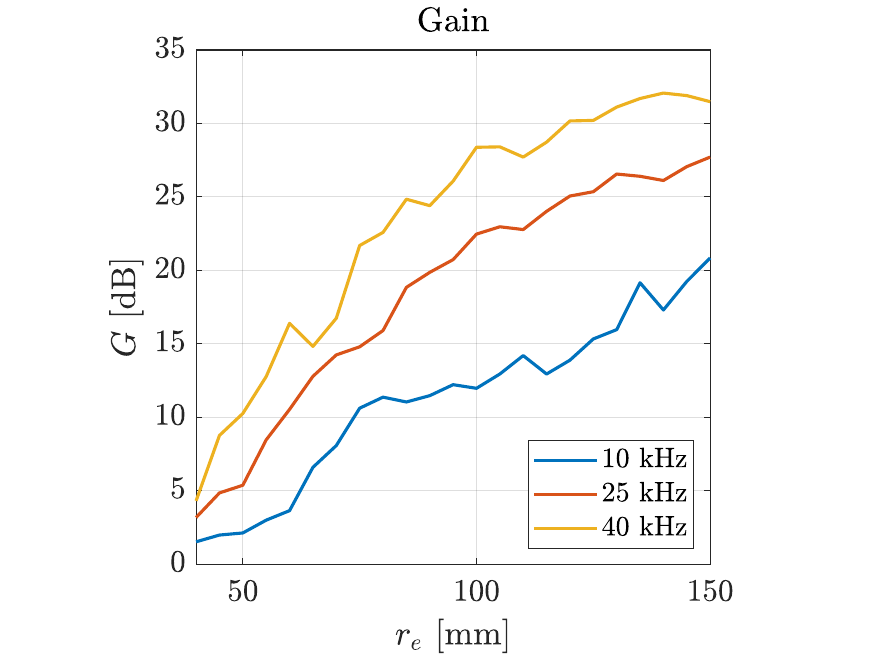}\label{fig:sizing costs-c}}
    \caption{Optimization results for the sizing procedure at the three frequencies. \protect\subref{fig:sizing costs-a} comparison of the performance achieved by the lens with smooth varying properties and by the discrete devices. The dashed line shows the expected quadratic trend of the four devices.
    \protect\subref{fig:sizing costs-b} values of the cost functional and \protect\subref{fig:sizing costs-c} gain with respect to the radius of the discrete lenses.
    }
    \label{fig:sizing costs}
\end{figure}

Figure~\ref{fig:sizing costs-a} shows the amplitude of the mean pressure in the focal region $|\langle p^i+p^s\rangle_{\Omega_f}|$, where 
\begin{equation}
    \langle\, p(\bb x)\,\rangle_{\Omega_f}\coloneqq |\Omega_f|^{-1}\int_{\Omega_f} p(\bb x)\,d\Omega \,\in\C.
\end{equation}
Note that the spatial average is an approximation of the measurement provided by a sensor placed in the middle of the device and it is proportional to the energy density of the signal in the focus.
% As expected, the gain notably increases with the external radius.
The curve indicates that the most substantial improvement occurs when the ratio of lens radius to wavelength, $r_e/\lambda$, approaches~$1$.

The energy conveyed in the focus by the device depends on its volume, hence the energy density is proportional to the ratio between the device and the focus volumes, i.e., $r_e^2/r_f^2$.
Thus, the gain behaves according to the dashed line shown in Figure~\ref{fig:sizing costs-a}.
% 
% An additional efficiency index is introduced to evaluate the performance of the device. Unlike a closed resonator, which can theoretically exhibit an infinite response, this device is open and lacks rigid boundaries. As a result, energy cannot be stored indefinitely, leading to leakage and reduced efficiency.
% We introduce the index
% \begin{equation}
%     \eta \coloneqq 1 - \frac{1}{|\langle p^i+p^s\rangle_{\Omega_f}|},
% \end{equation}
% such that $\eta=0$ in free-field and $\eta\to1$ for an ideal resonator.
% The blue dashed curve in Figure~\ref{fig:sizing costs-a} shows that the efficiency is higher that \qty{99}{\percent} when $r_e/\lambda>1$.

\subsection{Lattice discretization}
Assuming the long-wavelength homogenization applies \cite{laude2015phononic}, the metamaterial composing the lens is equivalent to an assembly of hexagonal cells having constant properties. Thus, the OCP is modified in order to account for the discrete nature of the device. The control space $\mathcal U$ is set to extend on a hexagonal grid such that the control functions are piecewise constant. Thus, the basis of $\mathcal U$ is $\{\psi_j(\bb x)\}_1^{N_c}$, where $N_c$ is the number of cells composing the lens and $\psi_j(\bb x)$ is the indicator function of the $j$\textsuperscript{th} cell. The control actions $(v(\bb x),u(\bb x))$ are fully described by the two vectors $\bb v,\bb u\in\R^{N_c}$ such that $v(\bb x) = \bs\psi(\bb x)^\top \bb v$, $u(\bb x) = \bs\psi(\bb x)^\top \bb u$.\\
Thus, the norm $\left|v,u\right|_{\mathcal U}$ is represented in the discrete context as $\bb v^\top (D+H)\bb v + \bb u^\top (D+H)\bb u$. Here, $D$ is a diagonal matrix with entries corresponding to the areas of the associated cells, and $H$ represents the Laplacian derived from the graph created by the cells' topology. In this graph, an edge exists between cells if they are adjacent. $H$ is defined as
\begin{equation}
H_{ij}=
\begin{dcases}
    |\Lambda_i| & \text{if } i=j 
    \\
    -1 & \text{if } j\in\Lambda_i
    \\
    0 & \text{otherwise}
\end{dcases},
\end{equation}
with $\Lambda_j$ the set of cells adjacent to the $j$\textsuperscript{th} cell and $|\Lambda_j|$ its cardinality---see \cite{cominelli2022design} for a more detailed derivation.\\
As a result, the cost function can be written as
\begin{equation}
    J(\bb v,\bb u,p^s) = \frac{\sigma}{2} (\bb{v}^\top \left(D+H\right) \bb{v} + \bb{u}^\top \left(D + H\right) \bb{u}) - \frac{1}{2}  G(p^s)
\end{equation}
% Note that due to the structure of $H$, $H+D$ is positive definite.
and the corresponding reduced gradients are
\begin{equation}
\begin{aligned}
% \nabla J_{v_j} &= \sigma \left(\bb v^\top (D+H)\right)_j  + \int_{\Omega_{c,j}}\Re\{a\nabla(p^s+p^i)\cdot\nabla\bar\lambda\}\,d\Omega = 0,
% \\
% \nabla J_{u_j} &= \sigma \left(\bb u^\top (D+H)\right)_j  - \int_{\Omega_{c,j}}\Re\{b\omega^2 (p^s+p^i)\,\bar\lambda\}\,d\Omega = 0,
% \text{ where $\Omega_{c,j}$ is the portion of $\Omega_c$ within the $j$\textsuperscript{th} cell.}
\nabla J_{\bb v} &= \sigma \bb v^\top (D+H)  + \int_{\Omega_c}\bs\psi^\top\,\Re\{ a\nabla(p^s+p^i)\cdot\nabla\bar\lambda\}\,d\Omega = \bs 0,
\\
\nabla J_{\bb u} &= \sigma \bb u^\top (D+H)  - \int_{\Omega_c}\bs\psi^\top\,\Re\{b\omega^2 (p^s+p^i)\,\bar\lambda\}\,d\Omega = \bs 0.
\end{aligned}
\end{equation}

In this setting, the OCP is solved for the three frequencies $f_L$, $f_M$, and $f_H$ and for several values of $r_e\in[\qty{40}{\milli\meter},\qty{150}{\milli\meter}]$. The results are displayed in Figure~\ref{fig:sizing costs} in terms of optimization cost and average gain in the focal region. As expected, at fixed $r_e$ the performance are higher for higher frequencies. By comparing the achieved results at fixed $r_e/\lambda$, the optimization at low frequency suffers the discrete nature of the metamaterial less than the high frequency one, as highlighted by Figure~\ref{fig:sizing costs-a}.\\

% It is not clear \textit{a priori} how the GRIN lens interacts with pressure waves since its lattice is composed by finite size cells.\\

\section{Multidirectional lens}
\label{sec:multi-directional}
In practical applications it is of great importance that a lens concentrates signals coming from several directions and waves in a defined frequency band, while all the lenses optimized in the previous section are intended for one frequency and one incident direction only.
We modify the OCP considering the interaction of the lens with respect to many monoharmonic signals, referred to as design waves coming from different directions. Their interaction with the lens is described by the same state equation~\eqref{eq:state_dyn}, where $\omega$, $\bb k$, and $p^i(\bb x)$ change according to the frequency of interest and the direction of propagation of the signal.\\
Given the time harmonic response of each design wave, a proper cost functional is needed that measures the performance of the lens over the whole design spectrum. The following multiobjective cost functional is introduced.
\begin{equation}\label{eq:multiobj cost}
    J_N(v,u,p^s_1,p^s_2,\dots,p^s_N) =  \frac{\sigma}{2} \left|(v,u)\right|^2_{\mathcal U} - \frac{1}{2}\prod_{n=1}^N G(p_n^s)^{1/N},
\end{equation}
where $p^s_n$ is the scattered pressure field given by the $n$\textsuperscript{th} design signal, characterized by harmonic frequency $\omega_n$, and wave vector $\bb k_n$, $k_n=\omega_n/c_0$, and $N$ is the number of design pressure waves considered.
Note that the geometric mean is taken over the $N$ intensities in order to promote a constant gain through the bandwidth.
In fact the geometric mean puts more weight on the smallest element.
In other words, the sensitivity computed with respect to an averaged element is higher for the smaller element. Thus the gradient based optimization updates the control keeping the gain uniform between different design waves.
Note that, if the algebraic mean is considered instead, the optimization would end up with a lens that gives priority to waves that require the smallest control effort.

The robust OCP, capable of optimizing the design spectrum, is
\begin{align}	\label{eq:cost functional robust}
	\min_{v,u,p^s_1,p^s_2,\dots,p^s_N} &\,J_N(v,u,p^s_1,p^s_2,\dots,p^s_N) \\
	s.t. &
	\begin{dcases} 	\label{eq:multiobj state_dyn}
		-\nabla \cdot (a\nabla p_n^s) - b \omega_n^2 p_n^s = f_n	& \text{in } \Omega
		\\
	    a\nabla p_n^s\cdot\mathbf{n} + \alpha_n p_n^s  =  0 	& \text{on } \partial\Omega
	\end{dcases}
\end{align}
where
\begin{equation}
    \begin{dcases}
		f_n = \omega_n^2(b-b_0)p_n^i +\nabla\cdot[(a-a_0)\nabla p_n^i]
		\\
		\alpha_n = a\Big(jk_n + \frac{1}{2r_a}\Big)
            \\
            p_n^i(\bb x) = e^{-j\bb k_n \cdot\bb x}
    \end{dcases}
\end{equation}

In order to address the $N$ constraints through the adjoint method, we define the augmented Lagrangian
\begin{equation}\label{eq:multiobj L}
	\LL_N \coloneqq J_N
	+ \sum_{n=1}^N \Re\Big\{\int_{\Omega} (\nabla \cdot (a\nabla p^s_n) + b \omega^2 p^s_n + f_n)\bar\lambda_n \, d\Omega 
    -\int_{\partial\Omega} (a\nabla p^s_n\cdot\bb n+\alpha_n p^s_n)\bar\lambda_n\,d\partial\Omega
    \Big\},
\end{equation}
where $\lambda_n\in \mathcal W^*$, for $n=1,2,\dots,N$ are the adjoint states.
Imposing the first order variation of $\LL_N$ with respect to the states $p^i_n$ to be null, the adjoint dynamics is derived
\begin{equation}\label{eq:multiobj adjoint}
    \begin{dcases}
        -\nabla \cdot (a\nabla \lambda_n) - b\omega^2_n \lambda_n = -(p^s_n+p^i_n)\, \chi_{\Omega_f}	G(p^s_n)^{-1} \prod_{m=1}^N G(p^s_m)^{1/N}	& \text{in } \Omega \\
        \nabla \lambda_n\cdot\mathbf{n} + \bar\alpha_n \lambda_n  =  0 		& \text{on } \partial\Omega
    \end{dcases}.
\end{equation}
Then, by computing the variation of $\LL_N$ with respect to the controls $v$ and $u$, the gradients of the reduced cost are obtained
\begin{equation}\label{eq:multiobj red_grad}
\begin{aligned}
    \nabla J_v &= \sigma (v+\Delta v) + \sum_{n=1}^N\Re\big\{a\nabla(p^s_n + p^i_n)\cdot\nabla \bar\lambda_n\big\}=0 \\
    \nabla J_u &= \sigma (u+\Delta u) -  \sum_{n=1}^N\Re\big\{b\omega^2_n (p^s_n + p^i_n) \,\bar\lambda_n\big\}=0,
\end{aligned}
\quad\text{in }\Omega_c,
\end{equation}
along with the natural boundary conditions $\nabla v\cdot\bb n = 0$ and $\nabla u\cdot\bb n =0$ on $\partial\Omega_c$.\\
Equations~\eqref{eq:multiobj state_dyn}, \eqref{eq:multiobj adjoint}, and \eqref{eq:multiobj red_grad} are the first order necessary conditions for the continuous setting.
Then, the discrete nature of the metamaterial is considered imposing the control space $\mathcal U$ to extend on a discrete hexagonal grid, similarly to the previous section.
% , and three lenses are optimized as follows.

We extend the optimizations performed previously considering three frequency bands centered on $f_L$, $f_M$, and $f_H$ and having an incident angle from \qty{-50}{\degree} to \qty[retain-explicit-plus]{+50}{\degree}.\\
At higher frequencies the ratio $l/\lambda$ is more critical, resulting in a coarser tuning of the lens properties. Thus, the width of the spectrum is kept constant and equal to \SI{2}{\kilo\hertz}, therefore that it corresponds to \qty{20}{\percent}, \qty{8}{\percent}, and \qty{5}{\percent} of the central frequency for the three bands, respectively. By doing so, similar performance are obtained by the three lenses.

A discrete number of frequencies and incidence angles must be selected for the robust optimization. Hence, six frequencies and six angles are selected such that they are equally spaced within the desired intervals. More precisely, the central frequency is modified by $\{\pm0.167, \pm0.5, \pm0.833\} \unit{\kilo\hertz}$ and the incident direction ranges between $\{\pm\qty{10}{\degree}, \pm\qty{30}{\degree}, \pm\qty{50}{\degree}\}$; thus, a total of $N=36$ design waves are enforced for each lens.

The size of each lens is chosen for each frequency as a reasonable compromise of the ratio $r_e/\lambda$, taking into account the curves shown in Figure~\ref{fig:sizing costs}. So, $r_e$ for the low, mid, and high frequencies is such that $r_{e,L}=r_{e,M}=\qty{74.3}{\milli\meter}$ and $r_{e,H}=\qty{55.7}{\milli\meter}$, namely $r_{e,L}/\lambda_L = 0.5$, $r_{e,M}/\lambda_M=1.25$, $r_{e,H}/\lambda_H= 1.5$.

\begin{figure}
    \centering
    \subfloat{\includegraphics[height=0.95\textwidth,trim=10 0 0 0,clip]{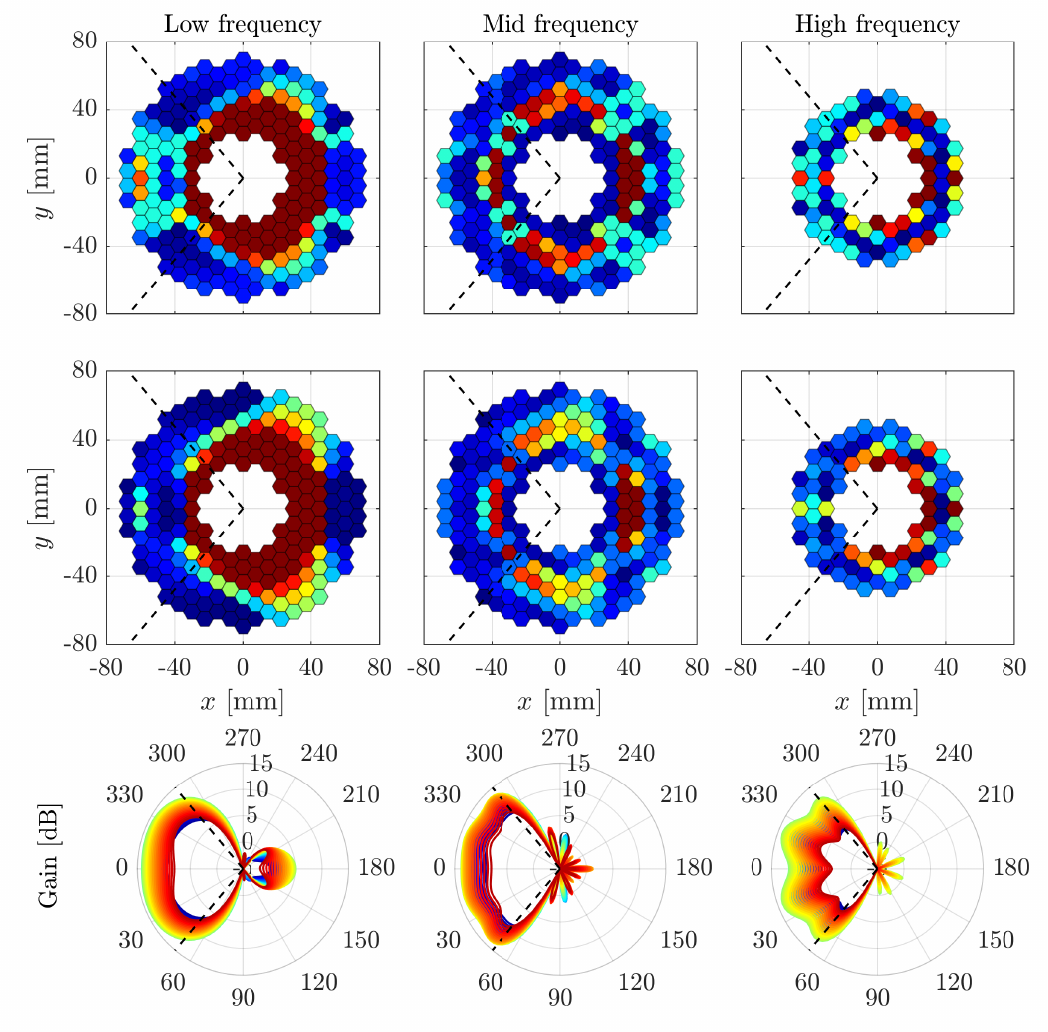}\label{fig:low freq lens-a}}
    \subfloat{\includegraphics[height=0.963\textwidth,trim=127 0 0 0,clip]{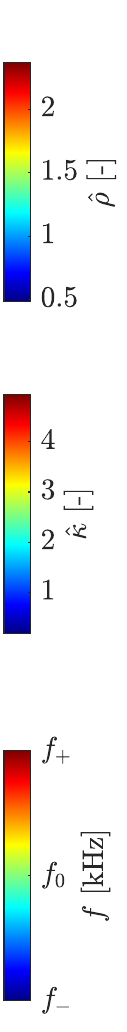}\label{fig:low freq lens-b}}
    \caption{Properties and performance of the optimized lenses. The dashed lines highlight the incidence acceptance cone of $\pm\qty{50}{\degree}$.
    % \seba{The three OCPs converged in $111$, $222$ and $333$ iterations, respectively.}
    For the sake of clarity, the polar plots do not show values below \qty{-5}{\decibel}.
    Note that, since the symmetry of the problem, the three lenses show a strong symmetry about the $x$ axis.
    }
    \label{fig:low-mid-high freq lenses}
\end{figure}

Finally, we set $\sigma=\num{1e-3}$ and solve the three OCPs. The properties and the performance of the three optimized lenses are shown in Figure~\ref{fig:low-mid-high freq lenses}.
The density and the bulk modulus of each cell are reported in color scale, while the amplitude of the mean pressure gain $|\langle\, p^i+p^s\,\rangle_{\Omega_f}|$ is displayed in three polar plots.\\
The polar cardioids are similar to circles within \qty[retain-explicit-plus]{+50}{\degree} and \qty{-50}{\degree}. This confirms the low directivity within this incident direction range. The frequency responses have a maximum in correspondence of the central frequencies and a quite flat behavior: from \qty{7.05}{\decibel} to \qty{14.35}{\decibel} for the low frequency,
from \qty{8.09}{\decibel} to \qty{14.01}{\decibel} for the mid frequency, and
from \qty{3.59}{\decibel} to \qty{14.12}{\decibel} for the high frequency. Figure~\ref{fig:TF} shows the same data in terms of amplitude and phase with respect to the frequency axis.

\subsection{Frequency response}
So far, the model adopted for the optimal design relies on two strong hypothesis: on the one hand, the time harmonic simplification limits the study to the steady state regime, where an infinite long signal withstands multiple reflections due to the lens structure;
on the other hand, the long-wavelength homogenization allows to consider each cell as an ideal fluid having an effective bulk modulus and mass density.
Due to the first assumption, the optimization tends to maximize the superposition of those reflections that take place in the focal region. However, such a procedure does not consider that a wave packet comprising a band of frequencies can experience distortion. 
Due to the second assumption, instead, an infinite, periodic lattice that interacts with a signal having a limited frequency band should be considered. Clearly, the optimized lenses do not satisfy such assumptions.\\
As a result, the performance of the lens should be verified.

In the following, we address the distortion looking at the gain and the phase of the transfer function measured between the input signal and the spatial mean of the signal within the focal region $\Omega_f$, that is:
\begin{equation}
    TF(\omega,\bb k) \coloneqq \frac{\langle\, p^i_{\omega, \bb k}(\bb x)+p^s_{\omega,\bb k}(\bb x)\,\rangle_{\Omega_f}}{p^i_{\omega, \bb k}(\bb x=\bs 0)},
\end{equation}
where $p^i_{\omega,\bb k}(\bb x)$, $p^s_{\omega,\bb k}(\bb x)$ are the complex valued incident and scattered pressure fields as a function of frequency and incidence direction.
% Knowing the transfer function, the amplified signal can be manipulated and the distortion recovered, nonetheless a contained distortion is preferred.
\\
Once again, the desired behavior of the lens is to concentrate the signal in its center, such that it experiences a constructive interference and the gain is constant between different frequencies and different directions of propagation, i.e.,  the amplitude $|TF(\omega,\bb k)|$ has to be as constant as possible.
When this is not the case, the so-called frequency response distortion occurs \cite{preis1982phase}. In this regard, the definition of the cost in \eqref{eq:multiobj cost} facilitates a flat response.

A further distortion comes from the relative phase shift between harmonics. By definition, there is no phase distortion if the zero-frequency intercept of the phase-frequency plot is $0$ or an integral multiple of $2\pi$ radians \cite{preis1982phase}.
In case of an odd multiple of $\pi$, a negative gain is identified and a simple signal inversion occurs.\\
In other words, an arbitrary time shift is chosen to interpret the received signal such that it has a small phase shift in the bandwidth of interest, according to the formula
\begin{equation}
    TF(\omega,\bb k) = \tilde{TF}(\omega,\bb k)\,e^{-j\omega\Delta t},
\end{equation}
where $\Delta t$ is the time delay and, in the phase-frequency plot of $TF$, it is interpreted as a slope, i.e.\ $\angle TF(\omega,\bb k) = \angle \tilde{TF}(\omega,\bb k) - \omega \Delta t$.

\begin{figure}
    \centering
    \raisebox{92pt}{\parbox{.03\textwidth}{\rotatebox{90}{
     $\angle \tilde{TF}$ [deg] \hspace{4pt}     
     $\angle TF$ [deg]   \hspace{4pt}     
     $|TF|$ [dB]
    }}}
    \subfloat{\includegraphics[height=0.44\textwidth,trim=10 10 20 0,clip]{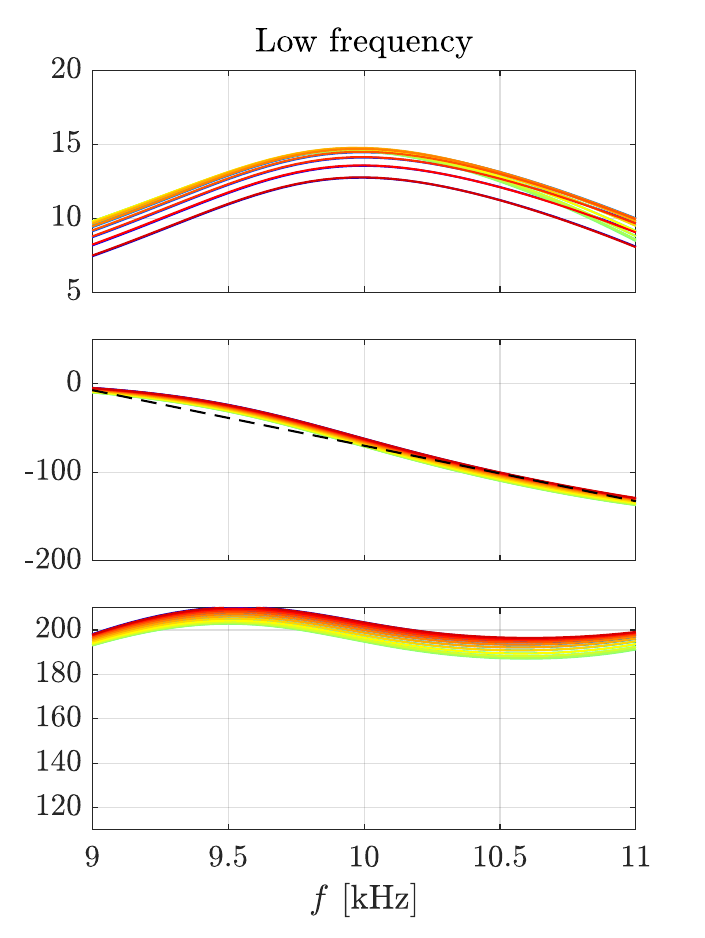}\label{fig:TF-a}}
    \subfloat{\includegraphics[height=0.44\textwidth,trim=35 10 20 0,clip]{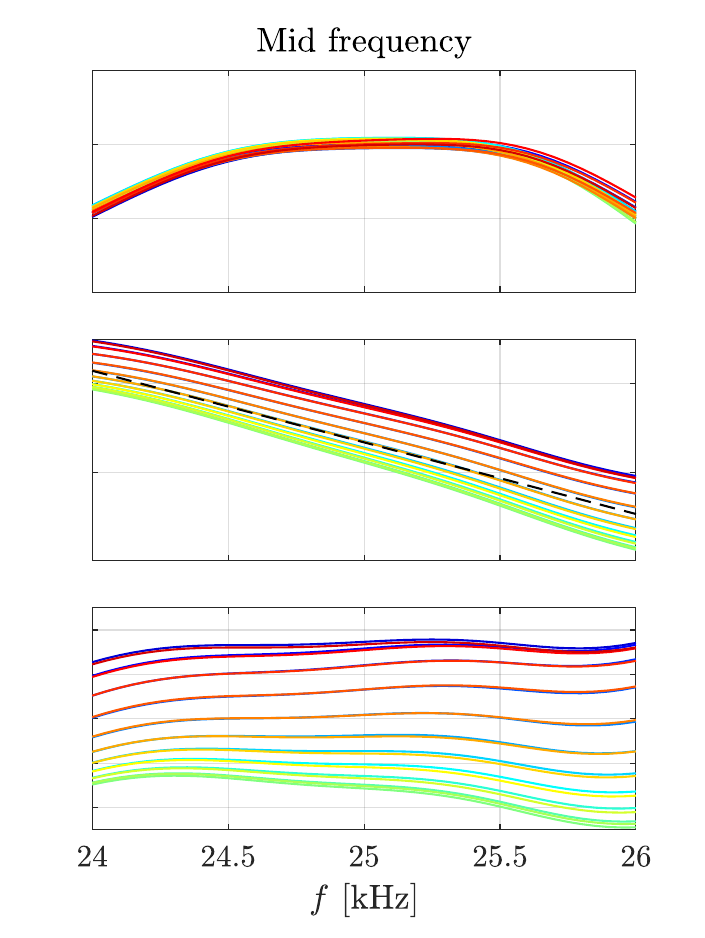}\label{fig:TF-b}}
    \subfloat{\includegraphics[height=0.44\textwidth,trim=35 10 20 0,clip]{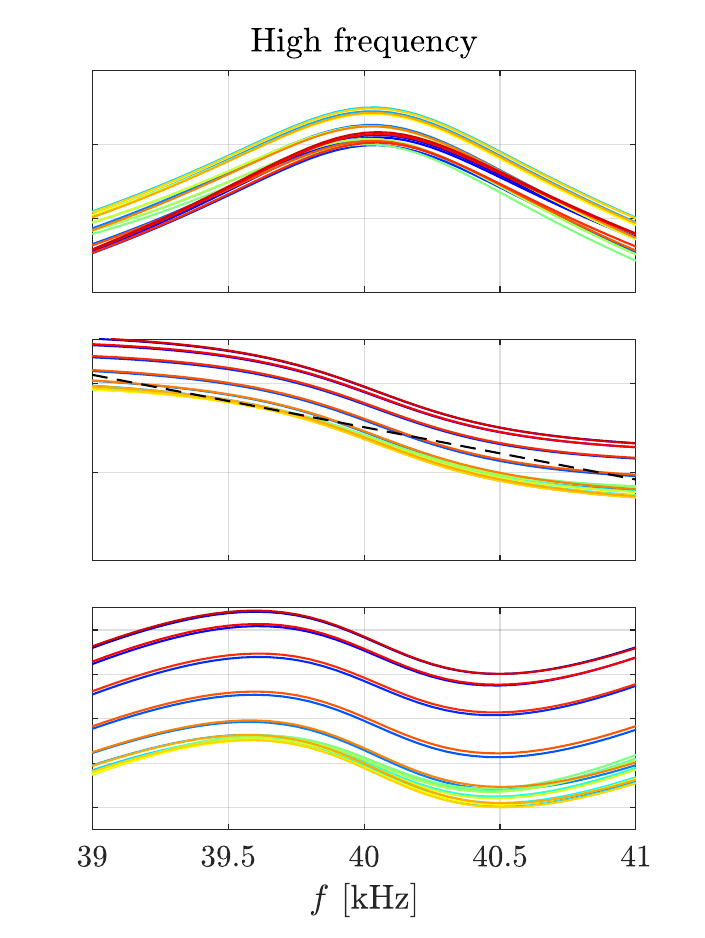}\label{fig:TF-c}}  
    \subfloat{  \includegraphics[height=0.44\textwidth,trim=265 10 10 0,clip]{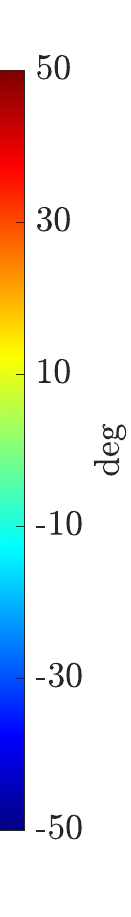}}
    \caption{From left to right. The first two lines shows amplitude and phase of the transfer function $TF$ of the three lenses; the bottom line shows the phase distortion.}
    \label{fig:TF}
\end{figure}

The first two rows of Figure~\ref{fig:TF} show amplitude and phase of the transfer functions granted by the three lenses in the respective frequency bands, the results obtained at different angles are shown in color scale.\\
All three gains have a broad peak and the signals with an incidence angle closer to \qty{0}{\degree} experience a slightly higher gain.
This condition would allow for a correction on the frequency response distortion even when the signal direction is unknown.\\
A relative phase shift can occur when signals propagate at different angles, potentially leading to destructive interference when multiple signals are combined. However, this is generally not a concern because the multi-path phenomenon inherently causes signals to arrive from various directions, often out of phase, even before any amplification occurs.

The phase distortion is addressed by looking at the phase-frequency plot. It shows an almost linear behavior and it is roughly approximated by the black line superimposed in the figures. Each straight line approximates the average of the phase through the different angles of incidence. The time delay $\Delta t$ is computed measuring the three slopes.

On the last column, the plots of $\angle\Tilde{TF}$ quantify the phase distortion of the signal.
The phase is about \qty{180}{\degree} for each lens. This means that the signal is inverted. On top of that, the low, the mid and the high frequency lenses show a phase distortion within $[\qty{15}{\degree},\qty{30}{\degree}]$, $[\qty{-70}{\degree},\qty{15}{\degree}$], and $[\qty{-60}{\degree},\qty{30}{\degree}]$ as a function of the incoming frequency, respectively. Such a sensible distortion shall be addressed by the electronics of the acquisition system.
% However, it is important to note that in all the three cases, a drastically smaller distortion occurs if only a narrower range of incidence direction is considered.

\subsection{Dynamics of the microstructure}

\begin{figure}
    \centering
    \subfloat{\includegraphics[width=0.33\textwidth,trim=0 0 0 0,clip]{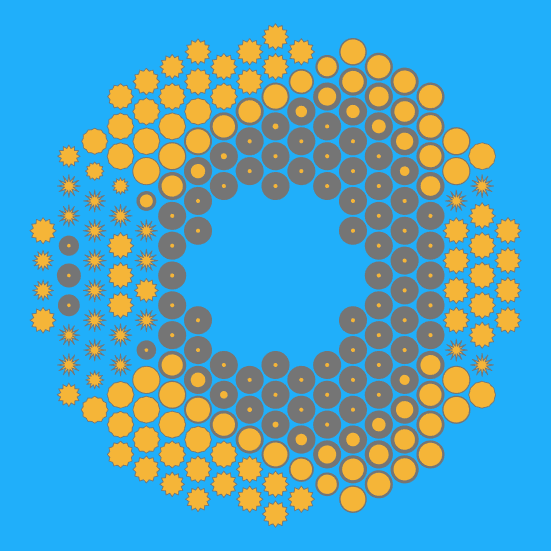}\label{fig:lente-a}}
    \,\subfloat{\includegraphics[width=0.33\textwidth,trim=0 0 0 0,clip]{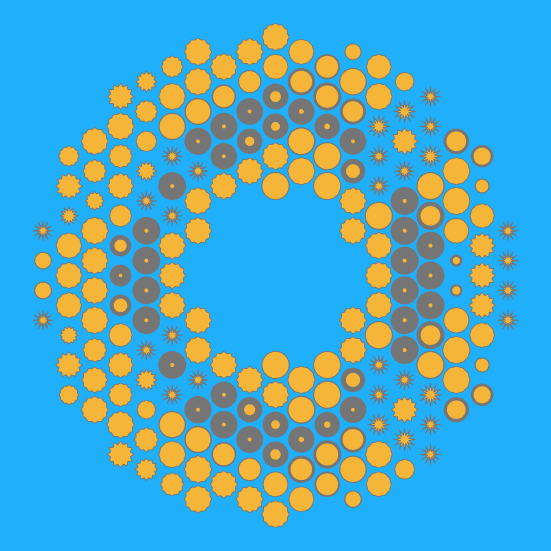}\label{fig:lente-b}}
    \,\subfloat{\includegraphics[width=0.33\textwidth,trim=0 0 0 0,clip]{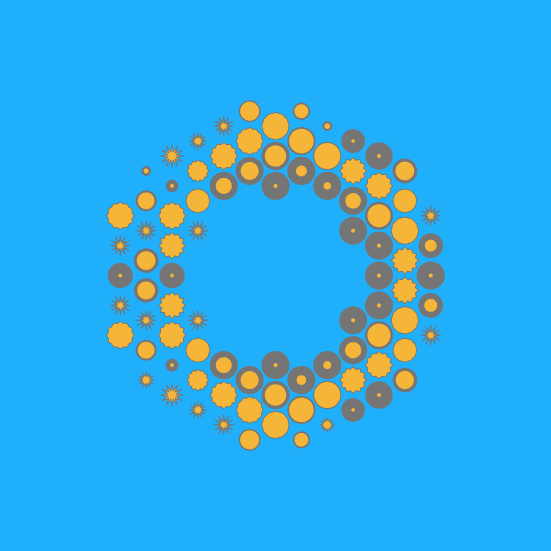}\label{fig:lente-c}}  
    \caption{From left to right, the microstructure of the low, mid and high frequency lenses.}
    \label{fig:lente}
\end{figure}

As anticipated, the homogenization applies to each cell if the lattice is periodic. However, due to the optimization each cell is different from the adjacent ones. Moreover, the dynamic behavior of the cells modifies the effective properties at high frequencies. Thus, a numerical simulation of the whole microstructure has to be carried out.\\
The first step is to convert the equivalent behavior of the cells in terms of $\rho$ and $\kappa$ into its geometric shape.
Thus, the attainable property space is explored letting the geometric parameters vary among 2000 different configurations. Then, a neural network (NN) capable of relating the properties $(\hat\rho,\hat\kappa)$ to either the parameters $(r_{in},r_{out})$ or $(p_{in},p_{out})$ has been trained. Since a smooth map $\R^2\to\R^2$ is needed, the corresponding NN is relatively simple to define and train, and it showed to be precise within a relative error on the properties less then \num{1e-4}. The resulting microstructures are shown in Figure~\ref{fig:lente}.

Once the geometry is entirely defined, the dynamic response has been obtained by a coupled acoustic-elastic FE simulation able to capture the local resonances of the lattice. Once more, the time harmonic approximation is exploited such that all occurring resonances are fully developed.
Please note that the signal is assumed to have a null incidence angle, i.e.\ $\bb k = (k,0)$.
% Hence, an investigation for all directions of incidence would be redundant.
Indeed, the result is mainly affected by the dynamics of the cells, and similar conclusions are drawn about the response to the other incidence angles.
% hence the lenses are excited only by signals coming from an angle of \qty{0}{\degree}. Note that the high number of different cells in a single lens is capable to describe how likely resonances occur in a generic device composed by these inclusions.
Such a numerical analysis is performed through the commercial software \Comsol{}.\\
The structure model has been provided with a small isotropic loss factor equal to \qty{1}{\percent} of its elastic modulus. This is expected to model the little structural damping.

\begin{figure}
    \centering
    \raisebox{80pt}{\parbox{.022\textwidth}{\rotatebox{90}{
     $\angle TF$ [deg] \hspace{15pt}   
     $|TF|$ [dB]
    }}}
    \subfloat{\includegraphics[height=0.38\textwidth,trim=10 10 20 0,clip]{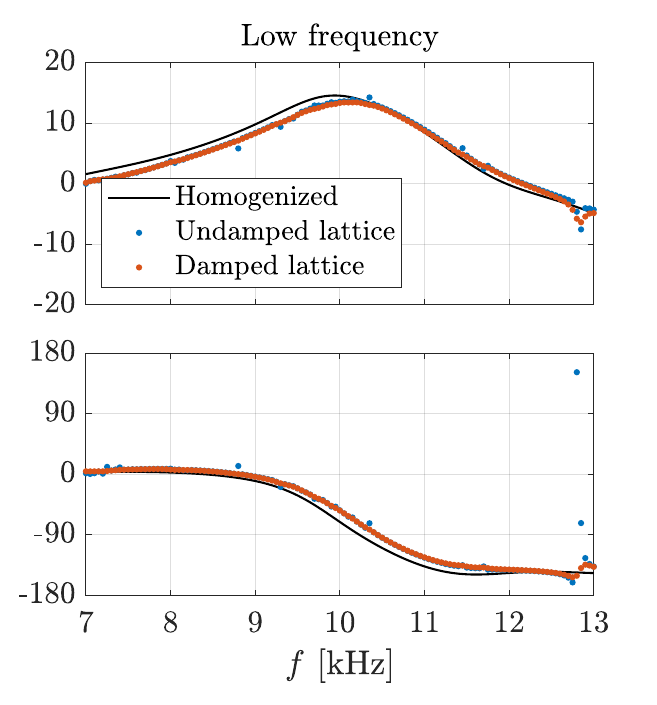}\label{figfig:microstructure validation-a}}
   \,\,\subfloat{\includegraphics[height=0.38\textwidth,trim=28 10 20 0,clip]{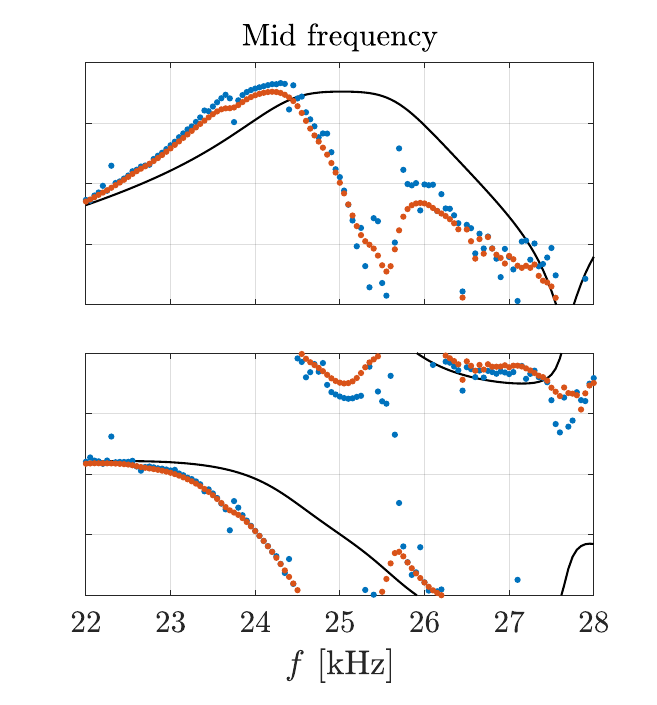}\label{fig:fig:microstructure validation-b}}
    \,\,\subfloat{\includegraphics[height=0.38\textwidth,trim=28 10 20 0,clip]{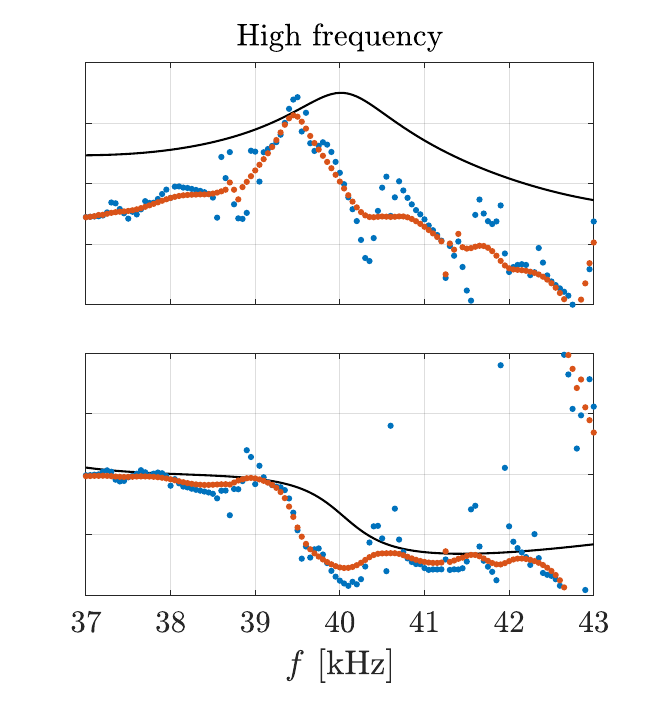}\label{fig:fig:microstructure validation-c}}
    \caption{Comparison between the transfer function of the homogenized lenses and their acoustic-elastic model. From left to right, the low, mid and high frequency lenses showing that the dynamics of the microstructure is progressively less negligible.
    The signal is assumed to have a null incidence angle, i.e.\ $\bb k = (k,0)$.}
    \label{fig:microstructure validation}
\end{figure}

The responses of the three lenses are compared with their ideal performance in Figure~\ref{fig:microstructure validation}.
The behavior of the low frequency lens shows a very good agreement in both undamped and damped approximations, and the transfer function is very close to that of the homogenized model.
The mid frequency lens, expected to work at about the theoretical limit of the cells' homogenization, shows a narrow peak shifted towards lower frequencies.
A similar and stronger effect occurs for the lens designed for high frequencies. Nonetheless, a sharp peak is evident at about \qty{39.5}{\kilo\hertz}.

Such an issue could be overcome if the high frequency behavior of the cell is accounted during the optimization.
However, this would require the so-called direct-method to evaluate the dispersion relation at a chosen frequency \cite{laude2015phononic},
% along with an optimization strategy that outlines the acoustic branch and evaluates the phase velocity.
which is out of the scope of this work.\\
Note that the performance of the mid and high frequency lenses are easily restored if a finer lattice is considered.

\section{Conclusions}
\label{sec:conclusion}

In marine environments, many technologies rely on manipulation of sound waves, because of their ability to propagate efficiently even for long distances.
In this work we investigate the performance of three GRIN lenses made of solid inclusions embedded in a water matrix.
Lenses working on a broadband spectrum and having a low directivity are targeted.\\
The lenses are composed by a graded metamaterial whose cells belong to a family of hallow scatteres.
Through efficient gradient-based optimization, we tune the homogeneous properties of each cell such that the grating behaves as a lens.
Given their acoustic properties, an efficient and precise NN performs the inverse design of each cell and the entire microstructure is accomplished.\\
The signal intensification provided by the lenses is discussed in terms of amplitude and phase distortion.

The scheme we propose is a valuable achievement in the control of acoustic waves because it is more powerful than those adopted in the literature, where the design procedure requires to discretize a smooth property variation \textit{a posteriori}, reducing the lens performance.
Conversely, a behavior close to reality is obtained at the end of this scheme.\\
We discuss the design of two-dimensional acoustic lenses, but the extension to a three-dimensional problem is straightforward from the optimization point of view, whilst it requires the design of a metamaterial suitable for the a 3D space.

The next important steps are to enrich the optimization scheme by introducing a cost on the phase distortion, to experimentally validate the lenses, and to improve an inverse design for high frequency homogenization.

% \appendix{\input{Sections/8_Appendix}}

\section*{Authors' contribute}
\textit{Sebastiano Cominelli}: Formal Analysis, Software, Validation, Visualization, Writing; \textit{Francesco Braghin}: Conceptualizaton, Supervision, Writing.

\section*{Competing interest}
The authors declare that they have no known competing financial interests or personal relationships that could have appeared to influence the work reported in this paper.

\section*{Funding}
The research presented in this paper was carried out within FRONTCOM (FRONTiers of Underwater COMmunications and Networking) project, co-funded by \textit{Segredifesa 5° Reparto – Innovazione tecnologica}.
% The research presented in this paper was supported by the FRONTCOM project—FRONTiers of Underwater COMmunications and Networking, funded under the SMARTCIG: ZF62DEC595 program by the Italian Navy.

\section*{Acknowledgements}
The authors thanks Davide E. Quadrelli, Matteo Tomasetto, and Alessandro Cominelli for the valuable discussions.

\bibliographystyle{RS}
\bibliography{main}

\begin{thebibliography}{99}

\bibitem{snelgrove2010discoveries}
Snelgrove PV. 2010 {\em Discoveries of the Census of Marine Life: making ocean
  life count}.
Cambridge University Press.

\bibitem{jiang2011electromagnetic}
Jiang S, Georgakopoulos S. 2011  Electromagnetic wave propagation into fresh
  water. {\em Journal of Electromagnetic Analysis and Applications}
  \textbf{2011}.

\bibitem{brown2003ray}
Brown MG, Colosi JA, Tomsovic S, Virovlyansky AL, Wolfson MA, Zaslavsky GM.
  2003  Ray dynamics in long-range deep ocean sound propagation. {\em The
  Journal of the Acoustical Society of America} \textbf{113}, 2533--2547.

\bibitem{zia2021state}
Zia MYI, Poncela J, Otero P. 2021  State-of-the-art underwater acoustic
  communication modems: Classifications, analyses and design challenges. {\em
  Wireless personal communications} \textbf{116}, 1325--1360.

\bibitem{furqan2019underwater}
Furqan~Ali M, Jayakody NK, D~Ponnimbaduge~Perera T, Srinivasan K, Sharma A,
  Krikidis I et~al.. 2019  Underwater communications: Recent advances. ETIC
  conference.

\bibitem{jiang2016convert}
Jiang X, Li Y, Liang B, Cheng Jc, Zhang L. 2016  Convert acoustic resonances to
  orbital angular momentum. {\em Physical review letters} \textbf{117}, 034301.

\bibitem{allam20203d}
Allam A, Sabra K, Erturk A. 2020  3D-printed gradient-index phononic crystal
  lens for underwater acoustic wave focusing. {\em Physical Review Applied}
  \textbf{13}, 064064.

\bibitem{lin2009gradient}
Lin SCS, Huang TJ, Sun JH, Wu TT. 2009  Gradient-index phononic crystals. {\em
  Physical Review B} \textbf{79}, 094302.

\bibitem{climente2010sound}
Climente A, Torrent D, S{\'a}nchez-Dehesa J. 2010  Sound focusing by gradient
  index sonic lenses. {\em Applied Physics Letters} \textbf{97}.

\bibitem{martin2010sonic}
Martin TP, Nicholas M, Orris GJ, Cai LW, Torrent D, S{\'a}nchez-Dehesa J. 2010
  Sonic gradient index lens for aqueous applications. {\em Applied Physics
  Letters} \textbf{97}.

\bibitem{zigoneanu2011design}
Zigoneanu L, Popa BI, Cummer SA. 2011  Design and measurements of a broadband
  two-dimensional acoustic lens. {\em Physical Review B} \textbf{84}, 024305.

\bibitem{li2014three}
Li Y, Yu G, Liang B, Zou X, Li G, Cheng S, Cheng J. 2014  Three-dimensional
  ultrathin planar lenses by acoustic metamaterials. {\em Scientific reports}
  \textbf{4}, 6830.

\bibitem{huang2021frequency}
Huang S, Peng L, Sun H, Wang Q, Zhao W, Wang S. 2021  Frequency response of an
  underwater acoustic focusing composite lens. {\em Applied Acoustics}
  \textbf{173}, 107692.

\bibitem{ma2022underwater}
Ma F, Zhang H, Du P, Wang C, Wu JH. 2022  An underwater planar lens for
  broadband acoustic concentrator. {\em Applied Physics Letters} \textbf{120}.

\bibitem{brambilla2024high}
Brambilla G, Cominelli S, Verbicaro M, Cazzulani G, Braghin F. 2024  High Bulk
  Modulus Pentamodes: the Three-Dimensional Metal Water. {\em arXiv e-prints}.

\bibitem{whewell1854cambridge}
Whewell W. 1854 {\em The Cambridge and Dublin Mathematical Journal} vol.~9.
E. Johnson.

\bibitem{luneburg1966mathematical}
Luneburg RK. 1966 {\em Mathematical theory of optics}.
Univ of California Press.

\bibitem{su2017broadband}
Su X, Norris AN, Cushing CW, Haberman MR, Wilson PS. 2017  Broadband focusing
  of underwater sound using a transparent pentamode lens. {\em The Journal of
  the Acoustical Society of America} \textbf{141}, 4408--4417.

\bibitem{pendry2000negative}
Pendry JB. 2000  Negative refraction makes a perfect lens. {\em Physical review
  letters} \textbf{85}, 3966.

\bibitem{rahm2008design}
Rahm M, Schurig D, Roberts DA, Cummer SA, Smith DR, Pendry JB. 2008  Design of
  electromagnetic cloaks and concentrators using form-invariant coordinate
  transformations of Maxwell’s equations. {\em Photonics and
  Nanostructures-fundamentals and Applications} \textbf{6}, 87--95.

\bibitem{gokhale2012special}
Gokhale NH, Cipolla JL, Norris AN. 2012  Special transformations for pentamode
  acoustic cloaking. {\em The Journal of the Acoustical Society of America}
  \textbf{132}, 2932--2941.

\bibitem{veselago1967electrodynamics}
Veselago VG. 1967  Electrodynamics of substances with simultaneously negative
  values of $\epsilon$ and $\mu$. {\em Usp. fiz. nauk} \textbf{92}, 517.

\bibitem{imamura2004negative}
Imamura K, Tamura S. 2004  Negative refraction of phonons and acoustic lensing
  effect of a crystalline slab. {\em Physical Review B} \textbf{70}, 174308.

\bibitem{cominelli2024isospectral}
Cominelli S, Vial B, Guenneau S, Craster RV. 2024  Isospectral open cavities
  and gratings. {\em Proceedings of the Royal Society A} \textbf{480},
  20230853.

\bibitem{antonakakis2014gratings}
Antonakakis T, Baida FI, Belkhir A, Cherednichenko K, Cooper S, Craster RV,
  Dem{\'e}sy G, Desanto J, Granet G, Gralak B et~al.. 2014  Gratings: Theory
  and numeric applications. .

\bibitem{shelby2001experimental}
Shelby RA, Smith DR, Schultz S. 2001  Experimental verification of a negative
  index of refraction. {\em Science} \textbf{292}, 77--79.

\bibitem{cominelli2022design}
Cominelli S, Quadrelli DE, Sinigaglia C, Braghin F. 2022  Design of arbitrarily
  shaped acoustic cloaks through partial differential equation-constrained
  optimization satisfying sonic-metamaterial design requirements. {\em
  Proceedings of the Royal Society A} \textbf{478}, 20210750.

\bibitem{cominelli2023optimal}
Cominelli S, Sinigaglia C, Quadrelli DE, Braghin F. 2023  Optimal strategies to
  steer and control water waves. {\em Ocean Engineering} \textbf{285}, 115346.

\bibitem{chen2021optimal}
Chen P, Haberman MR, Ghattas O. 2021  Optimal design of acoustic metamaterial
  cloaks under uncertainty. {\em Journal of Computational Physics}
  \textbf{431}, 110114.

\bibitem{bergmann1946wave}
Bergmann PG. 1946  The wave equation in a medium with a variable index of
  refraction. {\em The Journal of the Acoustical Society of America}
  \textbf{17}, 329--333.

\bibitem{banerjee2011introduction}
Banerjee B. 2011 {\em An introduction to metamaterials and waves in
  composites}.
Crc Press.

\bibitem{colton1998inverse}
Colton DL, Kress R, Kress R. 1998 {\em Inverse acoustic and electromagnetic
  scattering theory} vol.~93.
Springer.

\bibitem{schot1992eighty}
Schot SH. 1992  Eighty years of Sommerfeld's radiation condition. {\em Historia
  mathematica} \textbf{19}, 385--401.

\bibitem{bayliss1980radiation}
Bayliss A, Turkel E. 1980  Radiation boundary conditions for wave-like
  equations. {\em Communications on Pure and applied Mathematics} \textbf{33},
  707--725.

\bibitem{Comsolblog}
Frei W Using Perfectly Matched Layers and Scattering Boundary Conditions for
  Wave Electromagnetics Problems.
  \url{https://www.comsol.com/blogs/using-perfectly-matched-layers-and-scattering-boundary-conditions-for-wave-electromagnetics-problems/}.
Accessed: 2023-12-06.

\bibitem{Trol2010}
Tr{\"o}ltzsch F. 2010 {\em Optimal control of partial differential equations:
  theory, methods, and applications} vol. 112.
American Mathematical Soc.

\bibitem{laude2015phononic}
Laude V. 2015 {\em Phononic crystals: artificial crystals for sonic, acoustic,
  and elastic waves} vol.~26.
Walter de Gruyter GmbH \& Co KG.

\bibitem{benthowave}
 Omnidirectional Spherical Hydrophone.
  \url{https://www.benthowave.com/products/BII-7000omnihydrophones.html}.
Accessed: 2023-12-14.

\bibitem{nocedal2006numerical}
Nocedal J, Wright S. 2006 {\em Numerical optimization}.
Springer Science \& Business Media.

\bibitem{preis1982phase}
Preis D. 1982  Phase distortion and phase equalization in audio signal
  processing-a tutorial review. {\em Journal of the Audio Engineering Society}
  \textbf{30}, 774--794.

\end{thebibliography}

\end{document}